\documentclass[11pt]{amsart}
\usepackage{amssymb,hyperref,amsmath}

\usepackage{graphicx}
\newtheorem{theorem}{Theorem}%[section] %needed first because it defines the counter

\newtheorem{conjecture}[theorem]{Conjecture}
\newtheorem{corollary}[theorem]{Corollary}
\newtheorem{definition}[theorem]{Definition}
\newtheorem{example}[theorem]{Example}

\newtheorem{fact}[theorem]{Fact}
\newtheorem{lemma}[theorem]{Lemma}
\newtheorem{problem}[theorem]{Problem}
\newtheorem{proposition}[theorem]{Proposition}
\newtheorem{question}[theorem]{Question}
\newtheorem{remark}[theorem]{Remark}

\newcommand{\bcon}{\begin{conjecture}}
\newcommand{\econ}{\end{conjecture}}
\newcommand{\bcor}{\begin{corollary}}
\newcommand{\ecor}{\end{corollary}}
\newcommand{\bdf}{\begin{definition}}
\newcommand{\edf}{\end{definition}}
\newcommand{\beq}{\begin{equation}}
\newcommand{\eeq}{\end{equation}}
\newcommand{\bexa}{\begin{example}}
\newcommand{\eexa}{\end{example}}
\newcommand{\bfac}{\begin{fact}}
\newcommand{\efac}{\end{fact}}
\newcommand{\blem}{\begin{lemma}}
\newcommand{\elem}{\end{lemma}}
\newcommand{\bprb}{\begin{problem}}
\newcommand{\eprb}{\end{problem}}
\newcommand{\bpro}{\begin{proposition}}
\newcommand{\epro}{\end{proposition}}

\newcommand{\bque}{\begin{question}}
\newcommand{\eque}{\end{question}}
\newcommand{\brem}{\begin{remark}}
\newcommand{\erem}{\end{remark}}
\newcommand{\bthm}{\begin{theorem}}
\newcommand{\ethm}{\end{theorem}}
\newcommand{\bmat}{\begin{matrix}}
\newcommand{\emat}{\end{matrix}}

\newcommand{\bpr}{\begin{proof}}
\newcommand{\epr}{\end{proof}}

\newcommand{\comment}[1]{\,}

\newcommand{\R}{\mathbb R}

\title{A geometric invariant of virtual $n$-links}

\author{Blake K. Winter}
\thanks{The author wishes to thank the referee for pointing out an arithmetic error, and for many other very helpful suggestions, as well as Adam Sikora, for posing the question of whether the homotopy type of the Dehn space is a welded knot invariant.}
\address{Medaille College, Science, Mathematics and Technology\\
18 Agassiz Circle, Buffalo, NY, USA\\
bkwinter1985@gmail.com}

\begin{document}

\thispagestyle{empty}

%\begin{abstract}
%\end{abstract}
\begin{abstract}
For a virtual $n$-link $K$, we define a new virtual link $VD(K)$, which is invariant under virtual equivalence of $K$. The Dehn space of $VD(K)$, which we denote by $DD(K)$, therefore has a homotopy type which is an invariant of $K$. We show that the quandle and even the fundamental group of this space are able to detect the virtual trefoil. We also consider applications to higher-dimensional virtual links.
\end{abstract}

\pagestyle{myheadings}

\maketitle

%%%%%%%%%%%%%%%%%%%%%%%%%%%%%%%%%%%%%%%%%%%%%%%%%%%%%%%%%%%%%%%%
%

\section{Introduction}

Virtual links were introduced by Kauffman in \cite{Kauffman1}, as a generalization of classical knot theory. They were given a geometric interpretation by Kuperberg in \cite{Kuper}. Building on Kuperberg's ideas, virtual $n$-links were defined in \cite{BKW}. Here, we define some new invariants for virtual links, and demonstrate their power by using them to distinguish the virtual trefoil from the unknot, which cannot be done using the ordinary group or quandle. Although these knots can be distinguished by means of other invariants, such as biquandles, there are two advantages to the approach used herein: first, our approach to defining these invariants is purely geometric, whereas the biquandle is entirely combinatorial. Being geometric, these invariants are easy to generalize to higher dimensions. Second, biquandles are rather complicated algebraic objects. Our invariants are spaces up to homotopy type, with their ordinary quandles and fundamental groups. These are somewhat simpler to work with than biquandles.

The paper is organized as follows. We will briefly review the notion of a virtual link, then the notion of the \emph{Dehn space} of a virtual link. This is a space whose homotopy type is an invariant of the link. Next, we define the \emph{vertical double} of a virtual link, as well as the more general notion of a \emph{stack} for a virtual link. We then show that the fundamental group of the vertical double of the virtual trefoil distinguishes it from the unknot. Finally, we consider applications to higher-dimensional virtual links constructed by various spinning operations.

\section{Virtual \emph{n}-links}

A \emph{virtual $n$-link} is an embedding (smooth, PL, or TOP) $L:M^n \rightarrow F^{n+1} \times I$, where $M$ is an $n$-manifold, $F$ is an oriented $(n+1)$-manifold (possibly with boundary), and $I$ is a closed interval, with $L(M^n)\cap (F^{n+1} \times \partial I)$ empty. When possible without causing confusion, we simply refer to this virtual link as $L$, and we will abuse notation by neglecting to differentiate between the map and its image, when the distinction is clear. We consider two virtual $n$-links $L, L'$ equivalent by taking the equivalence relation generated by the following:
\begin{enumerate}
%	\item The images of $L, L'$ are identical in $F^{n+1} \times I$.
	\item Isotopy (smooth, PL, or TOP) of $L$ to $L'$.
	\item An orientation-preserving embedding $f \times id:F^{n+1} \times I \rightarrow F'^{n+1} \times I$  such that $f(L(M)) = L'(M)$.
\end{enumerate}

For $1$-links, a convenient set of diagrams may be created analogous to classical knot diagrams, but with one additional crossing type, the \emph{virtual crossing}. Virtual crossings are denoted by putting a small circle around the crossing, without cutting any of the arcs involved. \emph{Virtual equivalence} is generated by two types of moves: all classical Reidemeister moves are allowed on local regions of the diagram, and in addition, two diagrams are equivalent if one is obtained from the other by cutting out a sub-arc that only has virtual crossings and replacing it with another sub-arc that only has virtual crossings.

\begin{figure}
		\centering
			\includegraphics[scale=0.5]{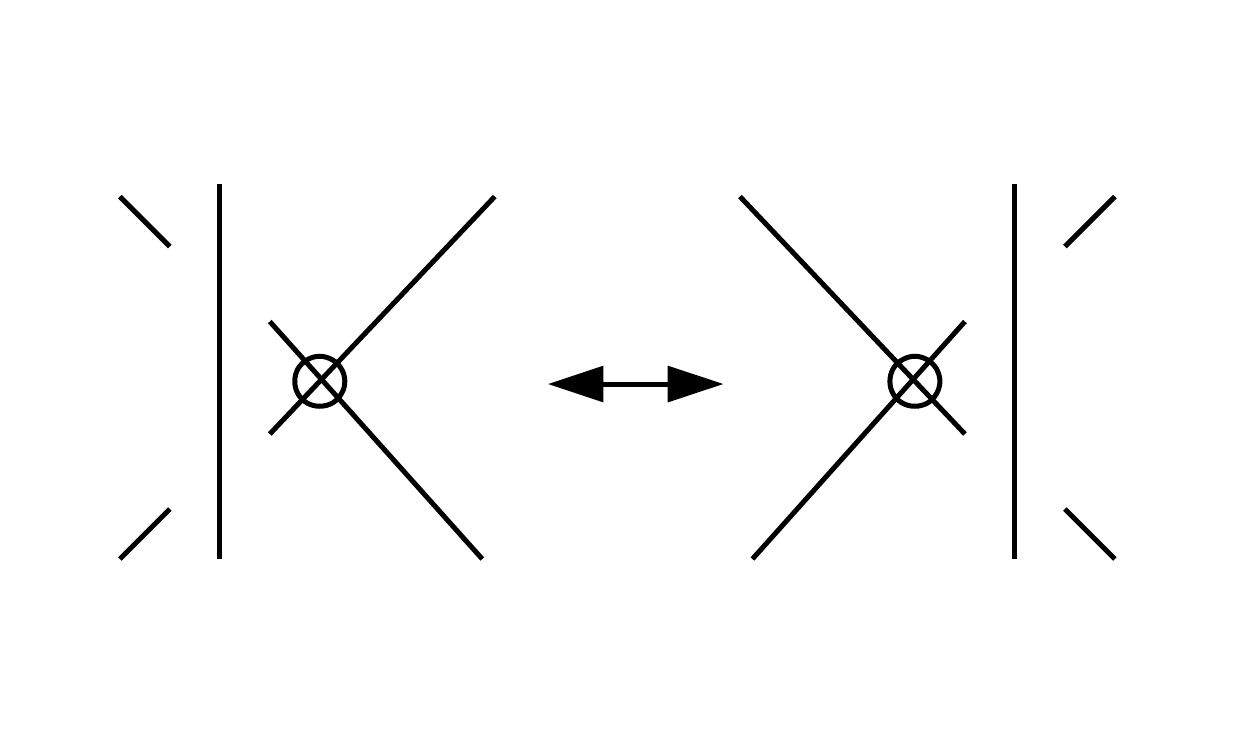}
		\caption{The \emph{welded forbidden move}.}
		\label{forbidden}
\end{figure}

There is another notion of equivalence for virtual $1$-link diagrams: \emph{welded equivalence}. This is generated by virtual equivalence together with the \emph{welded forbidden move}, shown in Fig. \ref{forbidden}. When considered up to the equivalence relation generated by virtual equivalence plus the welded forbidden move (or simply \emph{welded move}), we also refer to such links as \emph{welded links}. The welded move originated in the study of extensions of braid groups and the automorphisms of free quandles; see \cite{FRR}.

\section{The Dehn space of a virtual link}
An important invariant for virtual and welded links is the \emph{Dehn space}. The Dehn space is defined up to homotopy type. Let $L$ be a virtual $n$-link in $F \times I$. The Dehn space $D(L)$ is the space $D(L)=(F \times I - K)/F \times \{1\}$.

The fundamental group has a special subgroup, the \emph{peripheral subgroup}. Let $N(K)$ be the closure of a tubular neighborhood of $K$. Then the embedding of $\partial N(K)$ into $D(K)$ induces a homomorphism of fundamental groups. The image of this homomorphism is the peripheral subgroup. This subgroup has a special element (again, up to conjugation), the \emph{meridian}, which is the element that gets sent to $1$ under the Hurewicz homomorphism (note that a choice of generator for the first homology is equivalent to a choice of orientation for that component of the link).

\begin{theorem}
The homotopy type of $D(L)$ is invariant under virtual equivalence, under a homotopy that takes $\partial N(L)$ to $\partial N(L')$ and each meridian of $L$ to the corresponding meridian of $L'$.
\end{theorem}
\bpr
A detailed proof for all dimensions was given in \cite{BKW}, though for 1-links this result seems to be folklore. The result is easy to see when the images of two links are identical, as well as for isotopies of $L$, so it only remains to check what happens when we apply the third type of relation defining virtual equivalence. In this case we have an embedding $f \times id:F \times I \rightarrow F' \times I$  such that $f(L(M)) = L'(M)$. Then by sliding the points up, we can deformation retract $F'\times I$ so that it is contained in $f \times id (F \times I) \cup F' \times \{ 1\}$. Now after applying the quotient to obtain the Dehn space, we see that the homotopy types are the same.
\epr

The fundamental group of $D(L)$ is called the (\emph{virtual}) \emph{knot} (\emph{link}) \emph{group} of $L$, which we will denote by $G(L)$. We may also define the \emph{knot (link) quandle} of $L$, denoted $Q(L)$, following \cite{DJ, Mat}, using the Dehn space. Since the quandle is determined by the the fundamental group, the peripheral subgroup, and the meridian (see \cite{DJ, BKW}), and these are invariant under virtual equivalence, the quandle is a virtual $n$-link invariant.

\begin{figure}
		\centering
			\includegraphics[scale=0.5]{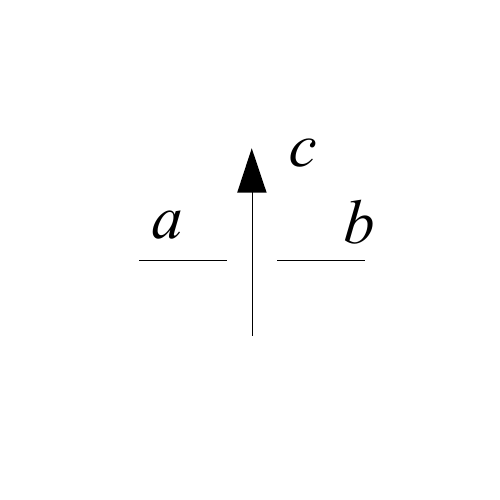}
		\caption{The Wirtinger relator for a double point as shown is $b=a^c$, where exponentiation denotes the quandle operation for quandles, and conjugation for groups. Herein we will use the convention that group multiplication in the fundamental group goes from left to right, so $b=cac^{-1}$. Note that this is the opposite of the convention used in some other publications. In higher dimensions, the generator $a$ is the one towards which the normal of the projection of $L$ points in the projection onto $F$. $c$ corresponds to the part of $L$ which is further from $F\times \{ 0 \}$.}
		\label{relator}
\end{figure}

To calculate the knot group or quandle, it is helpful to use the following procedure. Consider the projection $\pi$ of $F\times I$ onto $F$, and ensure that $L$ is in general position with respect to this projection. The projection is partitioned by its double point sets (those points which fail to be embedded). For each connected component of this partition, we get a generator. The double points give Wirtinger relators for the quandle or group, as shown in Fig. \ref{relator}. Details of this procedure, and a justification of it, may be found in \cite{BKW}; the proofs rely on a straightforward series of applications of the Van Kampen theorem, just as for classical links. The proof for classical links is well-known; see, for example, \cite{Rolfson}.

In the case of $1$-links, the homotopy type of $D(L)$ is invariant under the welded move. This means that it is not necessarily a strong invariant for virtual equivalence.

\begin{figure}
		\centering
			\includegraphics[scale=0.5]{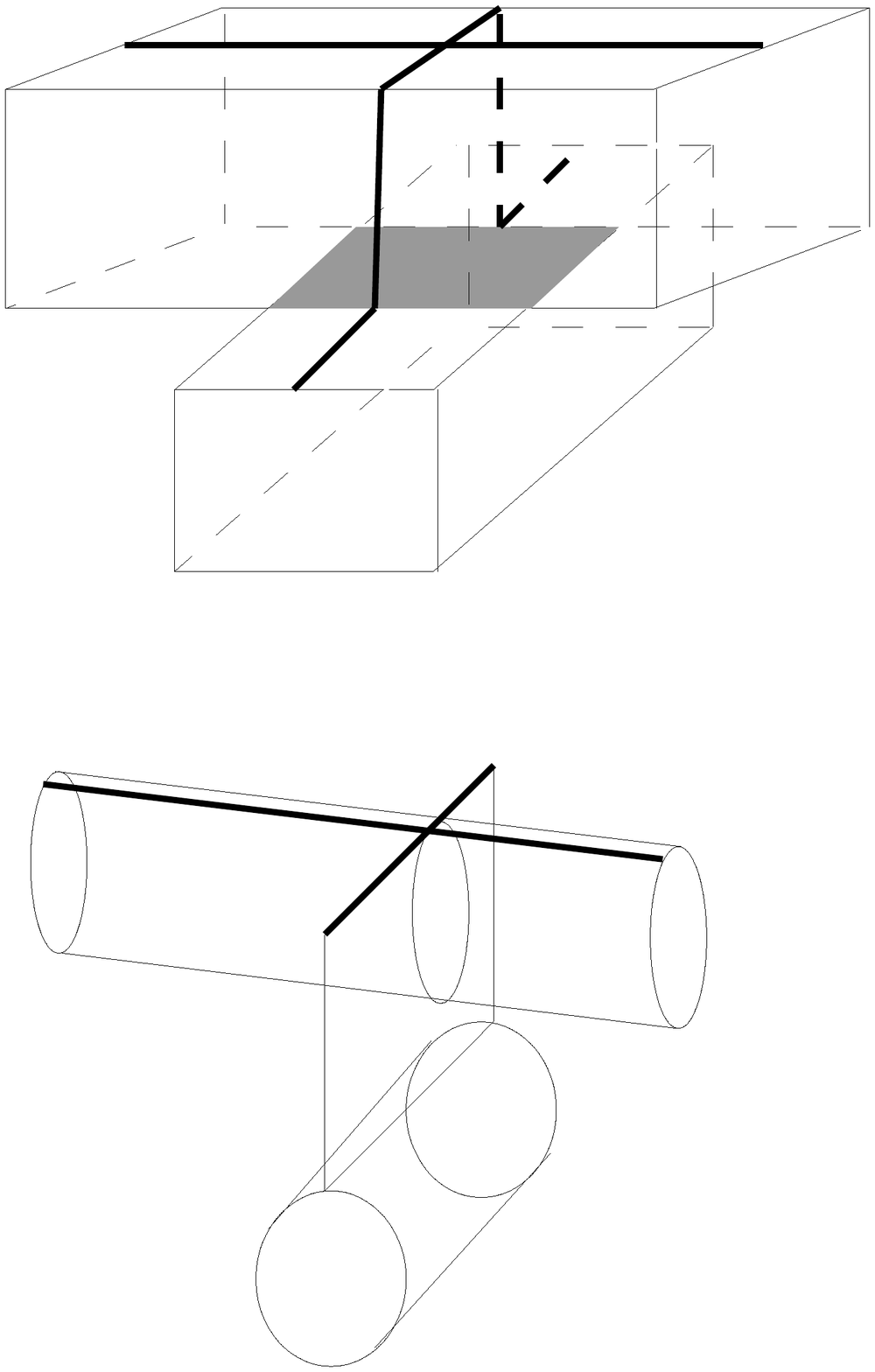}
		\caption{A deformation retraction of the Dehn space near a crossing. The darkened line shows the portion which is quotiented to a point.}
		\label{crossing}
\end{figure}

\begin{figure}
		\centering
			\includegraphics[scale=0.5]{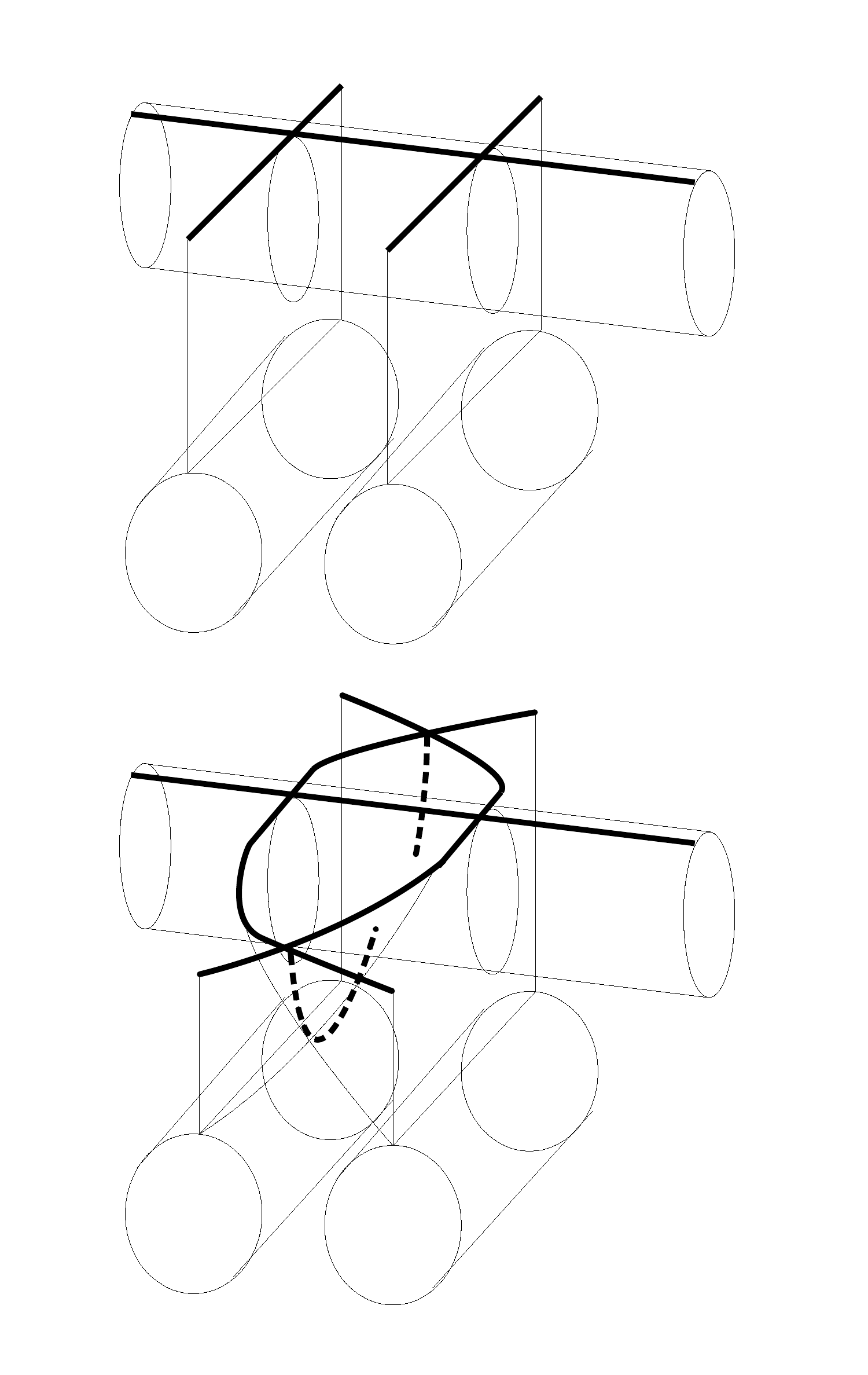}
		\caption{The change in the Dehn space due to a welded move. The darkened dotted line is an artifact of projecting into 3-space.}
		\label{crossingforbidden}
\end{figure}

\begin{theorem}
The homotopy type of $D(L)$ is unchanged for a virtual $1$-link under the welded move.
\end{theorem}
\bpr
The Dehn space of $L$ near a crossing is homotopy equivalent, by a deformation retraction, to a space as shown in Fig. \ref{crossing}. After a welded move, this changes the space near those two crossings as shown in Fig. \ref{crossingforbidden}. But this can be achieved by just sliding the sheets past one another, which is a homotopy equivalence.

\epr

\section{Stacks and the vertical double of a virtual link}

Let $L:M\rightarrow F\times [0, 1]$. By an isotopy (moving all points only along the direction of the interval) we may assume that the image of $L$ is contained in $F\times (\frac{1}{2}, 1]$. We will let $\pi$ be the natural projection to $F$ of $F\times [0, 1]$, and $t$ be the projection to $[0, 1]$. Let $M' = M \times \{0, 1\}$. We will define a new virtual link by $L':M'\rightarrow F\times [0, 1]$ as follows: $L'(x, 0)= (\pi(L(x)), t(L(x)))$,  $L'(x, 1)=(\pi(L(x)), 1-t(L(x)))$. 

Note that this simply means that we create a second copy of the link reflected over $F \times \{ \frac{1}{2} \}$. If we delete the original components of $L$ from $L'$, we would be left with the \emph{vertical mirror image} or vertical reflection of $L$. In terms of virtual link diagrams, this vertical reflection simply switches over and undercrossings.

We define the $L'$ virtual link to be the \emph{vertical double} of $L$, and denote it by $VD(L)$. We will use $DD(L)$ to refer to the Dehn space of $VD(L)$, as well as calling this the \emph{double Dehn space} of $L$. For the classical case, vertical doubles were introduced in \cite{IT}.

\begin{theorem}
If $L, K$ are virtually equivalent, then $VD(L)$ and $VD(K)$ are virtually equivalent, via a virtual isotopy sending the original components of $L$ to the original components of $K$.\label{VDequiv}
\end{theorem}
\bpr
Any isotopy of $L$ can be modified so that, throughout the isotopy, $L$ remains in $F\times (\frac{1}{2}, 1]$. Then the mirror image across $F\times \{ \frac{1}{2} \}$ gives an isotopy of the other half of $VD(L)$. On the other hand, if $f \times id: F \times I \rightarrow F' \times I$, then it is easy to see that $VD(F')$ is $f\times id(VD(L))$.
\epr

\begin{remark}
A similar argument shows that if the vertical mirror images of $L, L'$ are equivalent, then $L$ and $L'$ are equivalent, as well.
\end{remark}

\begin{corollary}
The homotopy type of $DD(L)$ (including the peripheral structures for each component), the link group of $VD(L)$, and the quandle of $VD(L)$ are invariants of $L$. In addition, if $L, L'$ are virtually isotopic, then there is an isomorphism of the link groups of $VD(L), VD(L')$ which maps the meridians and longitudes of each component in $L$ to a component in $L'$, considering $L, L'$ as sub-links of $VD(L), VD(L')$. \label{PC}
\end{corollary}

\begin{remark}
There is a choice to be made as to the orientation of the new components in the vertical double. We will generally use the convention that the new mirrored component is given the opposite orientation to its original. This agrees with the orientation induced by the higher-dimensional spin of which the vertical double turns out to be an equatorial slice; see Sec. \ref{SecSlice}
\end{remark}

\begin{remark}
Note that there is no reason to expect $VD(L)$ or $DD(L)$ to be invariant under welded equivalence for virtual $1$-links. As we will see below, they are not.
\end{remark}

We can generalize the above in the following manner. Let $L$ be a virtual link in $F\times I$. Let $L^*$ denote its vertical mirror. For every sequence $s$ of 0s and 1s beginning $s_1=1$, we can create the \emph{stack} $S(s, L)$. If $s$ has length $n$, then for each $i\in \{1, 2, 3, ..., n\}$, let $L_i=L$ if $s_i=1$, and let $L_i=L^*$ if $s_i=0$. We now glue the $F_i\times I$ together by gluing $F_i\times \{ 1\}$ to $F_{i-1}\times \{0 \}$, for $i=2, 3, ..., n$. By reparametrizing, we see that this is in fact a virtual link in $F\times I$. It has a copy of $L$ which is "over" all the other copies of $L$ and $L^*$. Observe that $VD(L)=S((1,0), L)$. Thus, the vertical double is a special case of the stacking construction. Note also that $S((0),L)$ would simply be the vertical mirror image of $L$.

\begin{theorem}
If $L$ and $L'$ are virtually equivalent, then $S(s, L)$ is virtually equivalent to $S(s, L')$ via a virtual equivalence sending $L_i$ to $L'_i$, for any sequence $s$ starting with 1.
\end{theorem}

The proof of this theorem is similar to the proof of Thm. \ref{VDequiv}.

\section{Application of the double Dehn invariants}
In this section, we will show that $DD(L)$ detects the nontriviality of the virtual trefoil. In fact, neither the quandle and fundamental group both distinguish the virtual trefoil from the trivial knot.

\begin{figure}
		\centering
			\includegraphics[scale=0.5]{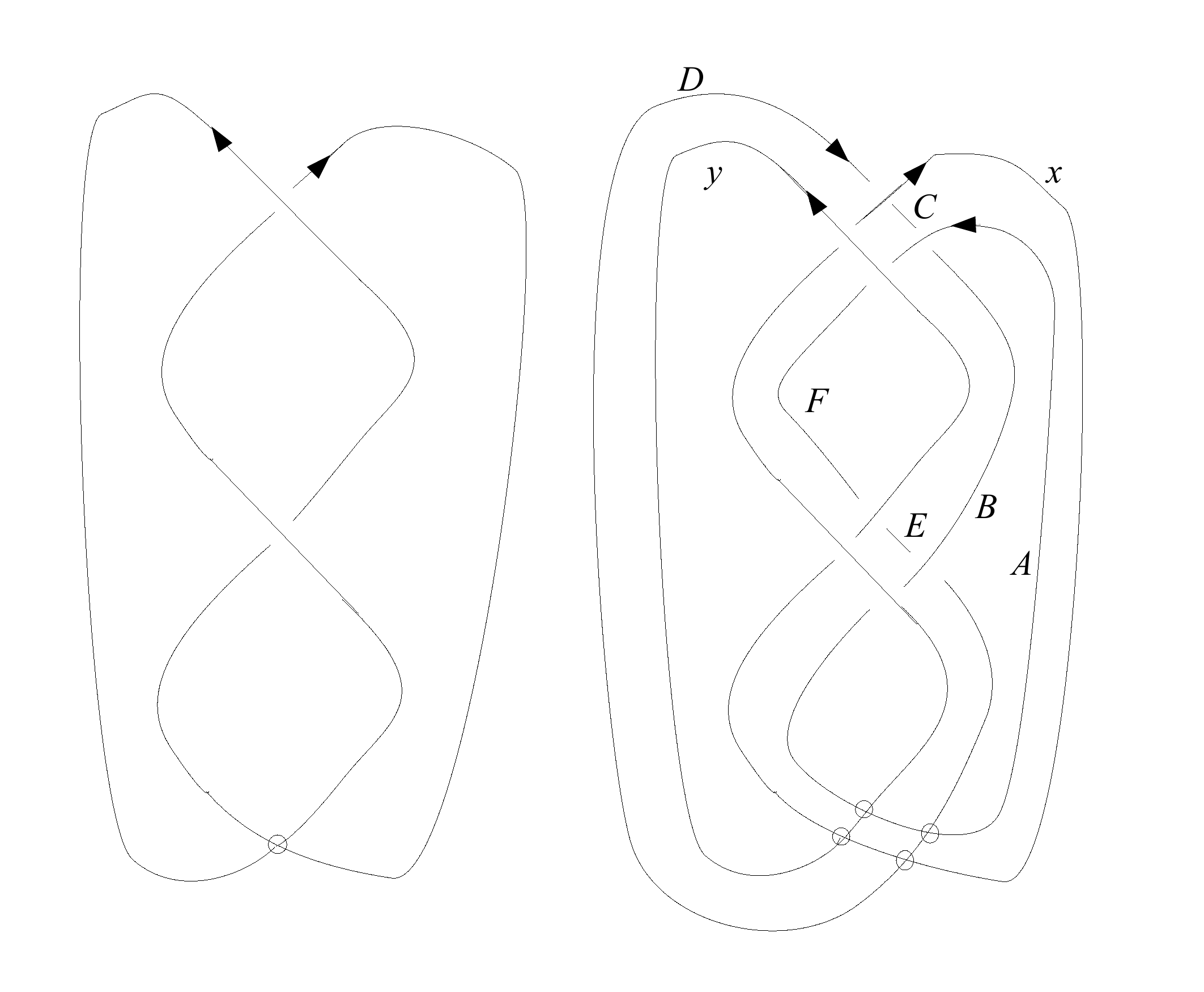}
		\caption{The virtual trefoil is shown on the left. The vertical double is shown on the right, with the arcs labeled to help calculate the quandle and group.}
		\label{VT}
\end{figure}

First, it is straightforward to see that $VD(U)$ for the unknot is the trivial 2-component unlink. Therefore, the quandle and group are the free quandle and group on two elements. We now turn to calculating the quandle for $VD(K)$ for $K$ the virtual trefoil. The virtual trefoil and its virtual double are shown in Fig. \ref{VT}.

From the labels in the figure, we are able to write down a presentation for the quandle. Using exponential notation, we have:

%	\[\left\langle x, y, A, B, C, D, E, F | x = y^y, y=x^y, C = D^x, C=B^A, A=F^y, E=F^y, E=D^B, B=A^y \right\rangle
%\]
\begin{align*}
\langle x, y, A, B, C, D, E, F | x = y^y, y=x^y, C = D^x, C=B^A, \\
A=F^y, E=F^y, E=D^B, B=A^y \rangle
\end{align*}
We may simplify this quandle presentation to obtain:
	\[\left\langle x, A, B, D, F | D^x=B^A, A=F^x, A=D^B, B=A^x \right\rangle
\]
	\[\left\langle x, A, B, D, | D^x=B^A, A=D^B, B=A^x \right\rangle
\]
	\[\left\langle x, A, D, | D^x=(A^x)^A, A=D^{(A^x)} \right\rangle
\]

This implies that $D=((A^x)^A)^{\overline{x}}$, using the overbar to indicate the inverse quandle operation. Therefore, we obtain the following simplified quandle presentation:

	\[\left\langle x, A | A = (((A^x)^A)^{\overline{x}})^{(A^x)} \right\rangle
\]

To see that this is not a free quandle, consider the free quandle on $x, A$. We construct a representation of the quandle into the quandle of $GL(2, \R)$ matrices with conjugation as the quandle operation. Note that this also induces a representation of the fundamental group into $GL(2, \R)$ as a group. The representation we choose will be:

	\[x\mapsto 
	\begin{bmatrix}
		x & y\\
		0 & z
	\end{bmatrix}, A\mapsto 
		\begin{bmatrix}
		a & b\\
		0 & c
	\end{bmatrix},
\]

Then a straightforward matrix arithmetic calculation shows that there is a nontrivial algebraic equation for $a, b, c, x, y, z$, in order for the relation $A = (((A^x)^A)^{\overline{x}})^{(A^x)}$ to hold. If we set $x=y=z=1, a=2, b=0, c=1$, for example, then we get the equation $0=-1$ in the upper right entries of the matrices. Therefore, for this choice of these variables, the relation does not hold in this homomorphic image of the free quandle, showing that our quandle (or group) is not free.

\begin{corollary}
The quandle and group of the double Dehn space detect the non-triviality of the virtual trefoil.
\end{corollary}

\begin{corollary}
The vertical double of a virtual knot, and double Dehn space homotopy type, is not a welded invariant.
\end{corollary}
\bpr
The virtual trefoil is descending and therefore trivial as a welded knot, but the vertical double distinguishes it from the unknot. See also \cite{SS2}.
\epr

\begin{remark}
The quandle and group of the vertical mirror image of the virtual trefoil are isomorphic to those of the unknot. Therefore, the vertical double carries more information about the virtual link than the vertical mirror image alone.
\end{remark}

The vertical double also provides us with a way to prove that a link is strictly virtual, that is, that it is not equivalent to any classical link.

\begin{theorem}
Let $L$ be a classical link with mirror image $L^*$. Then 
\[Q(VD(L))\cong Q(L)*Q(L^*).\]
\end{theorem}
\bpr
Since $L$ is classical, it is equivalent to a link in $D^{(n+1)}\times I$, where $D^{(n+1)}$ is an $(n+1)$-ball. When we form $VD(L)$, it is easy to see that there is a $(n+1)$-sphere which separates $L$ from the new components in $VD(L)$, and we may therefore perform an isotopy so that their projections to $D^{(n+1)}$ do not overlap. Then this result follows immediately.
\epr

For example, in the case studied above of the virtual trefoil, the quandles of the virtual trefoil and its mirror image are the trivial quandle, but the quandle of the vertical double is not the free product of two trivial quandles. Thus, the vertical double detects the non-classicality of the virtual trefoil.

\section{Equatorial links and two kinds of spinning} \label{SecSlice}

Let $L$ be an $n$-link in $F\times I$. We may assume that $L$ is in general position with respect to projection to $F$. We will form a new virtual $(n+1)$ link in the following manner. By an isotopy that does not change the projection to $F$, we may ensure that $L$ lies in $F\times (\frac{1}{2}, 1]$. Let $F'= F\times [-\frac{1}{2}, \frac{1}{2}]$. We obtain an $(n+1)$-link $L'$ by embedding $L\rightarrow F\times I\rightarrow F\times \{ 0\}\times I\subset F'\times I$, and then, in each copy of $[-\frac{1}{2}, \frac{1}{2}]\times I$, rotating the image of this embedding of $L$ about the point $(0, \frac{1}{2})$. If we think of $F\times[-\frac{1}{2}, \frac{1}{2}] \times I$ as a bundle over $F$, this simply rotates the image of $L$ inside each fiber. Where we used to have a crossing, we now have two copies of $S^1$, one inside the other. We will call $L'$ the \emph{vertical spin} of $L$. Then the vertical double of $L$ is obtained by considering the subspace of $F'\times I$ given by $(F\times \{ 0\})\times I$. In this sense, the vertical double is a kind of virtual equatorial link for the vertical spin of $L$ (for descriptions of classical equatorial links, see \cite{CKS}).

For $1$-links, this vertical spin is virtually equivalent to the Tube map defined by Satoh in \cite{SS}. However, the vertical double contains more information about $L$ than the vertical spin. This is because $\text{Tube}(L)$ is invariant under welded moves on $L$ (see \cite{SS, BKW}). It follows that the vertical spin of the virtual trefoil is virtually unknotted.

There is another notion of spinning virtual $n$-links given in \cite{BKW}. We define $\text{TSpun}(L)$ as follows. Given a virtual $n$-link $L$, $L:M\rightarrow F\times I$, let $M'=M\times S^1, F'=F\times S^1$. Then let $L':M'\rightarrow F'\times I$ be defined by $L'=L \times id_{S^1}$.

\begin{theorem}
$\text{TSpun}(VD(K))=VD(\text{TSpun}(K))$.
\end{theorem}
The proof of this theorem is an immediate consequence of the definitions. Note also that the TSpun operation preserves the knot group and quandle, \cite{BKW}. Therefore, the vertical double invariants distinguish the TSpun virtual trefoil from the vertically spun virtual trefoil, which is unknotted. Note that in \cite{BKW}, the TSpun virtual trefoil was shown to be nontrivially knotted using the biquandle.

The TSpun operation as given takes $1$-knots and creates knotted virtual tori. We can also form knotted spheres by a spinning procedure. Given a virtual $1$-link $L$ with underlying manifold $M$ in $F\times I$, suppose that $F$ has a non-empty boundary component. Choose an interval $B\subset \partial F$. Perform an isotopy on each component of $L$ so they each have exactly one interval lying in $B\times I$. Let $U$ denote the interior of $L^{-1}(B\times I)$. We now form a new space $F' = F\times S^1/R$, where $(x, t)R(x, t')$ for any $x \in B, t, t'\in S^1$. We form a virtual link $L'$ in $F'\times I$ by taking $L(M-U)\times S^1/R$. We will call $L'$ a \emph{sphere spin} of $L$, and denote this procedure by SSpun. This procedure produces a virtual 2-link whose components are all 2-spheres.

The calculations of the groups are facilitated by analyzing the double point crossings, as described in \cite{BKW, CKS}: the generators correspond to arcs of the projection of $L(M-U)$, while the relators are simply given by the Wirtinger relations from the crossings in the resulting virtual knotoids.

\begin{figure}
		\centering
			\includegraphics[scale=0.5]{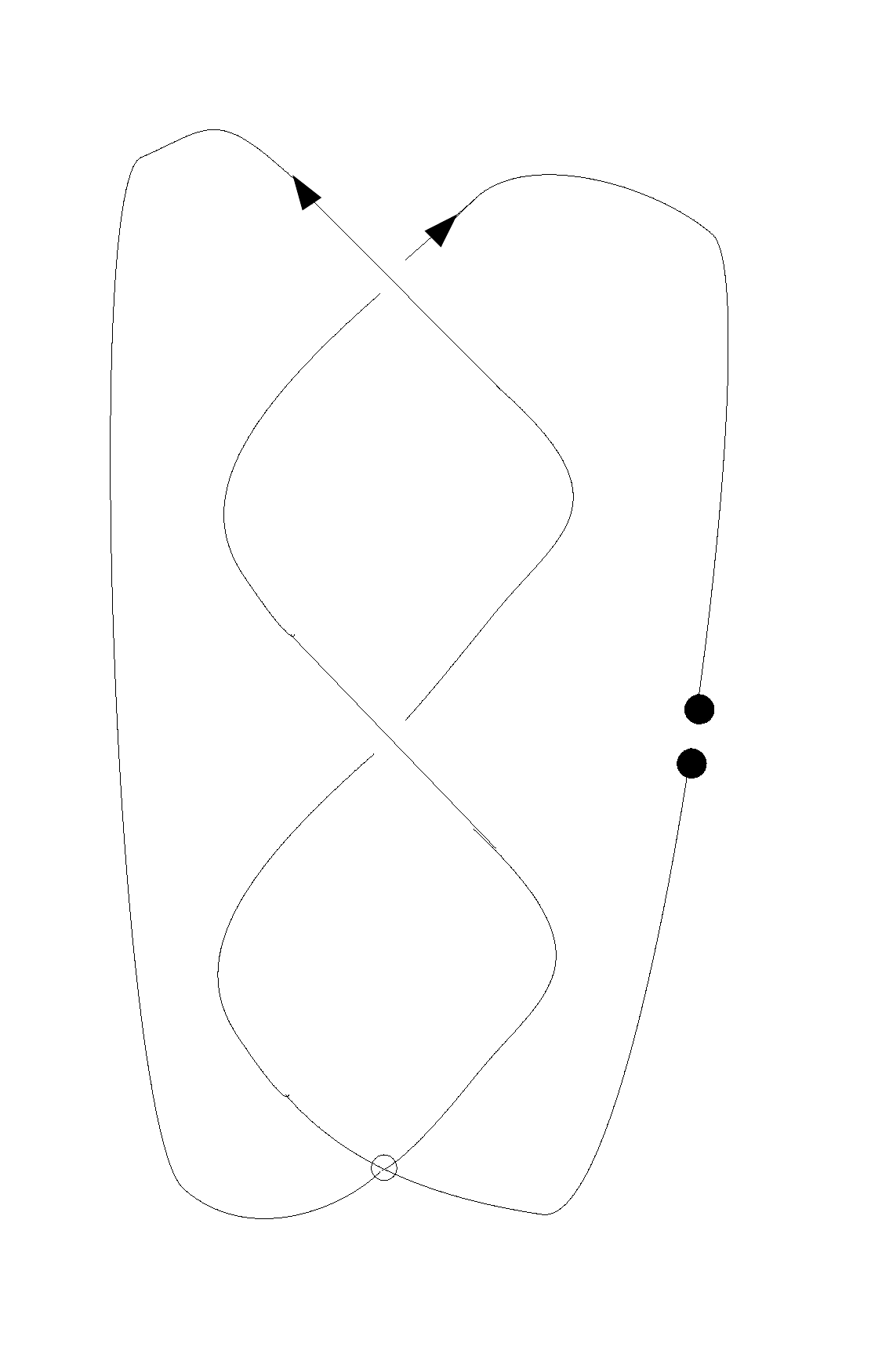}
		\caption{We can cut the virtual trefoil as shown and then spin it to form a sphere spun virtually knotted 2-sphere.}
		\label{trefspin}
\end{figure}

\begin{remark}
This operation depends upon the choice of where to cut the components of $L$ in general. For example, different choices for the virtual trefoil will produce virtually knotted 2-spheres with distinct fundamental groups: cutting the left arc gives the trefoil group, but cutting the right arc as in Fig. \ref{trefspin} yields the trivial group. Therefore, sphere spinning does not define a true operation on equivalence classes of virtual links. It requires us to choose a specific knotting and a point at which to cut it.
\end{remark}

\begin{remark}
In principle, for multicomponent links, we could choose different intervals on the boundary of $F$ to cut each component. This is achieved by allowing $B$ to be a disjoint union of intervals embedded in the boundary of $F$. However, for our purposes, we will not need this additional choice.
\end{remark}

Now consider performing the sphere spin operation on the virtual trefoil by cutting the arc between the marked points in Fig. \ref{trefspin}. Call this the \emph{sphere spun virtual trefoil}.

\begin{lemma}
The knot group of the sphere spun virtual trefoil is cyclic, and the quandle is the trivial quandle.
\end{lemma}

\begin{figure}
		\centering
			\includegraphics[scale=0.5]{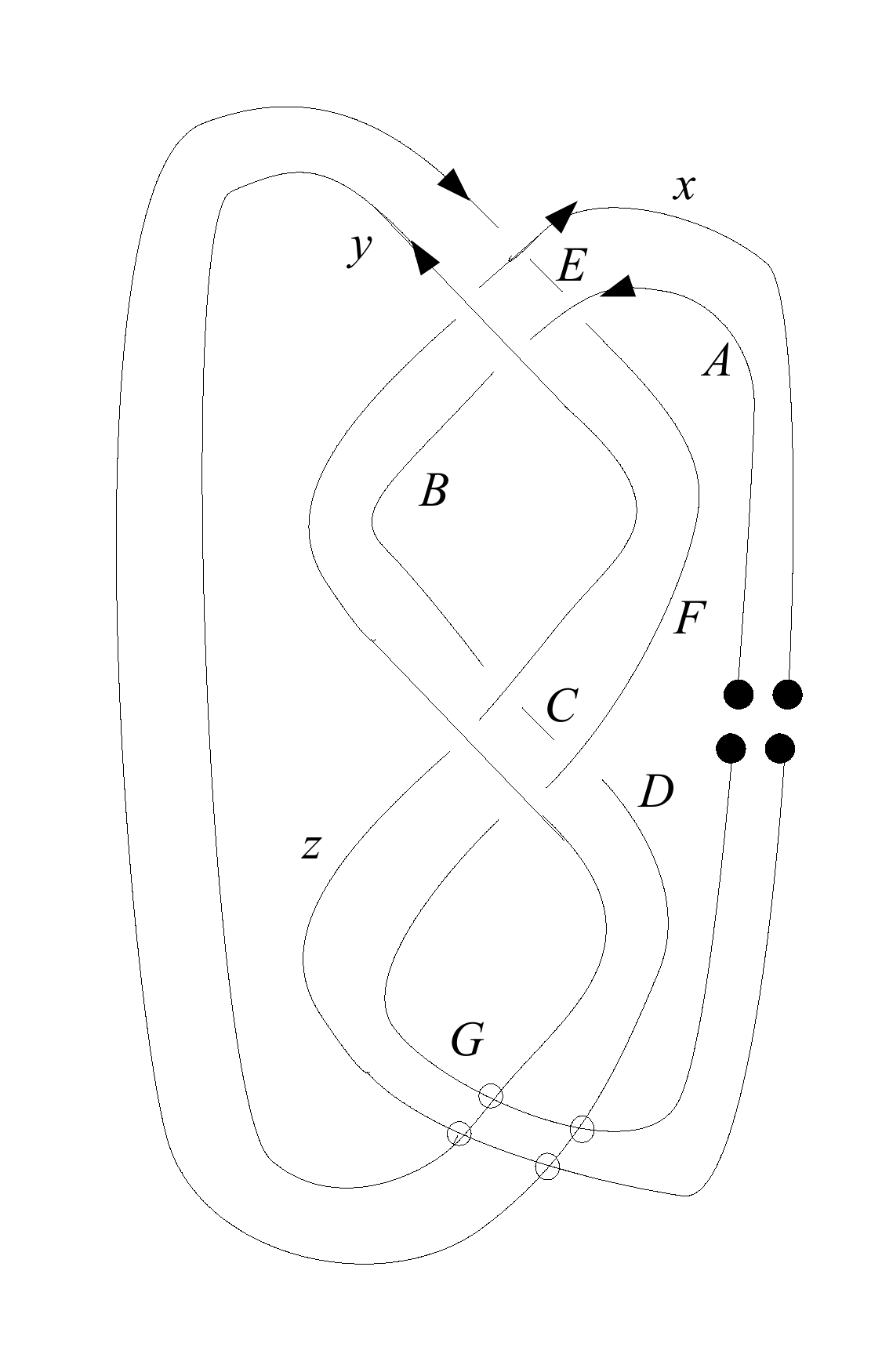}
		\caption{Cutting the vertical double of the virtual trefoil to spin it and obtain the vertical double of the sphere spun virtual trefoil.}
		\label{spundouble}
\end{figure}

\begin{theorem}
The vertical double of the sphere spun virtual trefoil is obtained by sphere spinning the vertical double of the virtual trefoil as shown in Fig. \ref{spundouble},
\end{theorem}

This theorem follows from the definitions of the vertical double and from the operation involved. The calculation of the fundamental group follows from analyzing the double point curves and the relations which they give.

We now show that the vertical double distinguishes the sphere spun virtual trefoil from the trivially knotted 2-sphere. Using the generators as shown in Fig. \ref{spundouble}, we obtain the following presentations for the quandle and group (for the group, interpret exponentiation as conjugation):

%	\[
%	\left\langle x, y, z, A, B, C, D, E, F, G | y=x^y, z=y^y, D=E^x, F=E^A, B=A^y, B=C^y, D=C^F, G=F^y \right\rangle
%\]

\begin{align*}
\langle x, y, & z, A, B, C, D, E, F, G | x=y^y, y=z^y, \\
	& E=D^x, E=F^A, A=B^y, C=B^y, C=D^F, F=G^y \rangle
\end{align*}

%	\begin{align*}
%	\left\langle x, y, & z, A, B, C, D, E, F, G | y=x^y, z=y^y, \\
%	& D=E^x, F=E^A, B=A^y, B=C^y, D=C^F, G=F^y \right\rangle
%\end{align*}

	\[\left\langle x, A, D, E, F, G | E=D^x, E=F^A, A=D^F, F=G^x \right\rangle
\]

	\[\left\langle x, A, D, F, G | D^x=F^A, A=D^F, F=G^x \right\rangle
\]
Now $D=(F^A)^{\overline{x}}$.
	\[\left\langle x, A, F, G| A=((F^A)^{\overline{x}})^F, F=G^x\right\rangle
\]
	\[\left\langle x, A, G| A=(((G^x)^A)^{\overline{x}})^{(G^x)}\right\rangle
\]

Note: at this point, if we set $G=A$, we obtain the virtual link group of the vertical double of the virtual trefoil, calculated previously. This is to be expected, since this would correspond to connecting the $A$ and $G$ arcs back up.

Now, observe that $x$ is the generator of the meridian of the spun virtual trefoil. If it were virtually trivial, then taking the quotient sending $x$ to $1$ would leave us with a cyclic group (see Cor. \ref{PC}). However, if we set $x=1$ in this presentation considered as a presentation of a group, we get the trefoil group. It follows that the vertical double distinguishes the sphere spun virtual trefoil from the trivially knotted 2-sphere.

This implies that we have a non-trivial virtually knotted $2$-sphere whose fundamental quandle is the trivial quandle, and whose fundamental group is cyclic. Since any classical knotted 2-sphere with cyclic fundamental group is trivial in the category TOP, \cite{CKS, Freedman}, this implies that, at least in the TOP category, we have constructed a virtual 2-knot which is strictly non-classical.

\begin{remark}
We could have seen the nontriviality of the sphere spun virtual trefoil by simply using the vertical mirror image and showing that it is nontrivial (it has a fundamental group isomorphic to that of the trefoil group). However, the group and quandle of the vertical double appear to contain \emph{more information} than this, since the generator $x$ is involved in a non-trivial relation with the generators $A$ and $F$.
\end{remark}

\begin{question}
The 1-knot known as the virtual trefoil has cyclic fundamental group and so does its vertical mirror image. Are there 2-knots which similarly have cyclic fundamental group and a vertical mirror image with cyclic fundamental group, but which are distinguished from a trivial knot by their vertical doubles?
\end{question}

For virtually knotted tori, the answer is easily seen to be yes. The virtual trefoil, spun to give a torus, gives an example. However, for virtually knotted spheres, we leave this as an open question.

\begin{remark}
If $K$ is a knot in $S^2 \times I$ (that is, a classical knot), then the sphere spinning construction simply gives Artin's spinning construction, discussed in \cite{CKS, SS}. In this case, the vertical double will be separable from $K$, and so will the vertical double of the sphere spun virtual 2-knot.
\end{remark}

\section{Conclusions and questions}
We observed that the quandle of a vertical double splits as the free product of the link and its vertical mirror when the link is classical. We may ask, therefore, whether this is a sufficient condition for classicality, as well. A somewhat stricter geometric question may also be posed:

\begin{question}
Suppose that $VD(L)$ is equivalent to the virtual link consisting of $L$ and its vertical mirror image modified so that their projections have no overlap. Does it follow that $L$ is classical?
\end{question}

The \emph{biquandle} of a virtual $n$-link, \cite{BQ1, biq, Carrell, BKW}, is another algebraic invariant which can distinguish the virtual trefoil from the unknot. However, the biquandle operations obey rather complicated relationships which may in some situations make working with the quandle or group of the vertical double of a virtual knot easier than working with the biquandle. It is therefore worth asking:

\begin{question}
If two virtual links are distinguished by their biquandles, does it follow that they can be distinguished by one of their stacking invariants? Or can a counterexample be constructed?
\end{question}

As an example of a virtual knot which is distinguished from the unknot by its biquandle, but not by the group of its vertical double, consider the Kishino knots studied in \cite{NelVo}. A somewhat lengthy computation shows that the link group of the vertical doubles of these knots are free groups on two generators. But it may be possible to distinguish them using other, more complicated stacking invariants. It would also be interesting to ask whether there is a limit to new information that can be obtained by taking stacking invariants involving longer sequences.


\begin{thebibliography}{9999}
%
%%%%%%%%%%%%%%%%%%%%%%%%%%%%%%%%%%%%%%%%%%%%%%%%%%%%%%%%%%%%%%%%%%%
\bibitem{Carrell}{T. Carrell, \emph{The surface biquandle}, Thesis, Pomona College (2009) available online http://www.math.washington.edu/~tcarrell/pomona-thesis.pdf.}
\bibitem{CKS} J.S. Carter, S. Kamada, and M. Saito, \emph{Surfaces in 4-space},
Springer 2004.
\bibitem{biq} R. Fenn, M. Jordan, L. Kauffman, \emph{Biracks, biquandles and virtual knots}, Topol. Appl. \textbf{145} (1-3) (2004).

\bibitem{FRR}{R. Fenn, R. Rimanyi, C. Rourke, \emph{The braid permutation group}, Topology \textbf{36} 1, (1997), 123-135.}
\bibitem{Freedman} M. H. Freedman, \emph{The disk theorem for four-dimensional manifolds}, Proc. Internat. Congr. Math. (Warsaw, Poland, 1983), PWN, Warsaw 
(1984) 647-663. 
\bibitem{IT} K. Ishikawa, K. Tanaka, \emph{Quandle colorings vs. biquandle colorings}, preprint: https://arxiv.org/abs/1912.12917.
\bibitem{DJ} D. Joyce, \emph{A Classifying Invariant of Knots, the Knot Quandle},
Journal of Pure and Applied Algebra \textbf{23} 37-65 (1982).
\bibitem{Kauffman1} L.H. Kauffman, \emph{Virtual knot theory}, European J. Comb. \textbf{20} (1999), 663-690.
\bibitem{BQ1} L. Kauffman, D. Radford, \emph{Bi-oriented quantum algebras, and a generalized alexander polynomial for virtual links}, Contemp. Math. \textbf{318} (2002), 113-140, available online http://homepages.math.uic.edu/\~kauffman/GenAlex.pdf.
\bibitem{Kuper} G.  Kuperberg, \emph{What is a virtual link?}, Algebr. Geom. Topol. \textbf{3} (2003), 587-591.

%\bibitem{Lee} J. Lee, \emph{Introduction to smooth manifolds}, Springer (2006).
%\bibitem{Liv}{C. Livingston, \emph{Stably irreducible surfaces in $S^{4}$}, Pacific Journal of Mathematics, Vol. 116, No. 1 (1985).}
\bibitem{Mat} S. Matveev, S.V. \emph{Distributivnye grupoidy v teorii uzlov}, Mat.
Sbornik \textbf{119} 1 (1982) 78-88 (in Russian). English Version: \emph{Distributive groupoids in Knot Theory}, Mat.
Sbornik \textbf{47} 73-83 (1984).
\bibitem{NelVo} Nelson, S., Vo, J., \emph{Matrices and finite biquandles}, Homology, Homotopy Appl. 8, (2006), 51–73.
\bibitem{Rolfson} D. Rolfson, \emph{Knots and links}, AMS Chelsea Publishing, 2003.
\bibitem{SS2} S. Satoh, \emph{Crossing changes, delta moves and sharp moves on welded knots}, Rocky Mountain J. Math.
    \textbf{48} 3 (2018), 967-879.
\bibitem{SS} S. Satoh, \emph{Virtual knot presentation of ribbon torus-knots},
J. Knot Theory Ramifications \textbf{9} (2000), 531-542.
%\bibitem{Naka} Y. Nakanishi, \emph{On ribbon knots, II}, Kobe. J. Math 3 (1986) 77-85.
%\bibitem{NakaNaka} Y. Nakanishi, Y. Nakagawa, \emph{On ribbon knots}, Math. Sem. Notes Kobe Univ. 5 (1982) 423-430.
\bibitem{BKW} B.K. Winter, \emph{Virtual links in arbitrary dimensions}, J. Knot Theory Ramifications \textbf{24} 14, 2015.

\end{thebibliography}
\end{document}